\documentclass{amsart}
\usepackage{amsmath,amsthm,amssymb,IMjournal,xypic}
\numberwithin{equation}{section}\theoremstyle{plain}
\newtheorem{The}{Theorem}[section]
\newtheorem{Def}[The]{Definition}
\newtheorem{Cor}[The]{Corollary}
\newtheorem{Lem}[The]{Lemma}
\newtheorem{Prop}[The]{Proposition}

\begin{document}

\title{Ding projective complexes with respect to a semidualizing module$^{\dagger}$}

\author{\|Yanhong |Quan|, \|Renyu |Zhao$^{\ast}$|, \|Chunxia |Zhang|, Lanzhou Gansu}


\dedicatory{Cordially dedicated to ...}

\abstract
Let $R$ be a commutative ring
and $C$ a semidualizing $R$-module. In this article, we introduce and
investigate the notion of $D_{C}$-projective complexes. We first prove that a complex $X$ is
$D_{C}$-projective if and only if each degree of $X$ is a $D_{C}$-projective module and $\mathcal{H}\text{om}(X,H)$ is exact for any $C$-flat
complex $H$. As immediate consequences of this result, some properties of $D_{C}$-projective complexes are given. Secondly, we investigate a kind of stability of $D_{C}$-projective complexes by showing that an iteration of the procedure used to define the $D_C$-projective complexes yields exactly the $D_C$-projective complexes. Finally, We introduce and characterize the notion of $D_C$-projective dimension of complexes.
\endabstract



\thanks
\hskip -0.4 true cm $^{\dagger}$ Supported by National Natural Science Foundation of China (11361052,11401475,11401476).
\newline $^{\ast}$Corresponding author.\endthanks

\section{Introduction and Preliminaries}\label{sec1}

In recent years, Gorenstein homological algebra has been developed to an advanced level, see for example \cite{ch,cfh,ejbook, h} and literatures list in them. It's main idea is to replace projective (resp. injective, flat) modules by Gorenstein projective (resp. Gorenstein injective, Gorenstein flat) modules. These modules were introduced by Enochs et al. \cite{ej,ejt} as generalizations and dauls
of finitely generated modules of G-dimension zero over a two-sided Noetherian ring in the sense of
Auslander and Bridger \cite{ab}. At the same time, These concepts have been extended in several directions. One of generalizations is Gorenstein modules with respect to a semidualizing module. As a generalization of dualizing modules and free modules of rank 1, Foxby \cite{f}, Golod \cite{g} and Vasconcelos \cite{v} independently initiated the study of semidualizing modules (under different names) over a commutative Noetherian ring. In particular, Golod \cite{g} used these to define $G$-dimension with respect a semidulizing module for finitely generated modules. Motivated by Enochs and Jenda's ideas in \cite{ej,ejt}, Holm and J${\o}$rgensen \cite{hj} focused on Gorenstein projective (resp. Gorenstein injective, Gorenstein flat) modules with respect to a semidualizing module $C$ over a commutative Noetherian ring, which were called $C$-Gorenstein projective (resp. $C$-Gorenstein injective, $C$-Gorenstein flat) modules. White \cite{w} extended the notions of semidualizing modules and Holm and J${\o}$rgensen's $C$-Gorenstein projective modules to commutative non-Noetherian rings, and she called $C$-Gorenstein projective modules as $G_C$-projective modules where $C$ is a semidualizing module. Many classical results about the Gorenstein projectivity of
modules were generalized in \cite{w}. Dually, $G_C$-injective modules were studied in \cite{yl}. Holm and White \cite{hw} further extended the definition of semidualizing modules to a pair of arbitrary associative rings, and many results on semidualizing modules over commutative Noetherian rings were generalized to this more general setting. In \cite{lhx}, the $G_C$-projective modules and the $G_C$-projective dimension of modules over general rings were investigated. In a different direction, Enochs and Garc\'{i}a Rozas \cite{eg,egr}
introduced Gorenstein projective (resp. injective, flat) complexes, and proved that over Gorenstein rings, these complexes are actually the
complexes of Gorenstein projective (resp. injective, flat) modules. Yang \cite{y} further proved that the Gorenstein projective (resp. injective) versions of the above results are true over arbitrary rings, and the Gorenstein flat version holds over coherent rings. Yang and Liang \cite{yl} introduced Gorenstein projective (resp. injective) complexes with respect to a semidualizing module $C$ over commutative rings, and proved that these
complexes are actually the complexes of $G_C$-projective (resp. injective) modules.

On the other hand, in \cite{dlm,dm}, Ding et al. considered two special cases of the Gorenstein
projective and Gorenstein injective modules, which they called strongly Gorenstein
flat and Gorenstein FP-injective modules, respectively. These two classes of modules
over coherent rings possess many nice properties analogous to Gorenstein projective
and Gorenstein injective modules over Noetherian rings. For the reason
that these modules were introduced and studied by Ding and his co-authors, Gillespie renamed strongly Gorenstein flat as Ding
projective, and Gorenstein FP-injective as Ding injective. In \cite{zwl}, Zhang, Wang and Liu introduced and studied Ding projective (resp. injective) modules with respect to a semidulaizing module over commutative rings. Ding projective (resp. injective) complexes were
investigated by Yang, Liu and Liang \cite{yll}, among others, they proved that over any
ring $R$, a complex $X$ is Ding projective (resp. injective) if and only if each $X^i$ is a Ding projective (resp. injective)
module for all $i\in \mathbb{Z}$ and $\text{$\mathcal{H}$om}(X,F)$ (resp. $\text{$\mathcal{H}$om}(J,X)$)is exact for any flat complex $F$ (resp. any FP-injective complex $J$).

Motivated by the above works, in this paper, we introduce and investigate Ding
projective (resp. injective) complexes with respect to a semidulaizing module. We only deal with Ding
projective complexes with respect to a semidulaizing module, Ding injective version can be given dually.

Next we shall recall some notions and definitions which we need in the later
sections. In order to make things less technical, throughout this article, by a ring $R$, we always mean a commutative ring with identity, all modules are unitary $R$-modules. We use ${\rm Ch}(R)$ to denote the category of complexes of $R$-modules.

\vskip.25cm{\bf~1.1} A complex
\vspace*{-2mm}$$\xymatrix{\cdots \ar[r]^{} &X^{n+1} \ar[r]^{\delta^{n+1}} &  X^{n}  \ar[r]^{\delta^{n}}& X^{n-1} \ar[r]^{}&\cdots}\vspace*{-2mm}$$
will be denoted by $(X,\delta)$ or simply $X$. The $n$th cycle (resp. boundary, homology) of
$X$ is denoted by $Z_n(X)$ (resp. $B_n(X)$, $H_n(X)$). Given an $R$-module $M$, we will denote by $\overline{M}$ the complex
\vspace*{-2mm}$$\xymatrix{\cdots \ar[r]^{} &0 \ar[r]^{} & M \ar[r]^{id}& M\ar[r]^{}&0 \ar[r]^{}&\cdots}\vspace*{-2mm}$$
with $M$ in the 1 and 0th degrees. Given an $X\in {\rm Ch}(R)$ and an integer $m$,
$X[m]$ denotes the complex such that $X[m]^n=X^{n-m}$ and whose boundary operators
are $(-1)^m\delta^{n-m}$. Given $X,Y\in {\rm Ch}(R)$, We use $\text{Hom}(X,Y)$ to present the group of all morphisms
from $X$ to $Y$, and Ext$^{i}(X,Y)$ denote the right-derived functors of Hom.  We let $\text{$\mathcal{H}$om}(X,Y)$ denote the complex
with \vspace*{-2mm}$$\text{$\mathcal{H}$om}(X,Y)^n=\prod_{t\in \mathbb{Z}}\text{Hom}(X^t,Y^{n+t}),\vspace*{-2mm}$$
and with differential given by
\vspace*{-2mm}$$\delta^n\left((f^t)_{t\in \mathbb{Z}}\right)=\left(\delta_{Y}^{n+t}f^t-{(-1)^n} f^{t-1}\delta_{X}^{t}\right)_{t\in \mathbb{Z}}.\vspace*{-2mm}$$

\vskip.25cm\textbf{~1.2} (\cite[1.8]{w}) An $R$-module $C$ is called semidualizing if

(1)~$C$ admits a degreewise finite projective resolution.

(2)~the natural homothety map $\xymatrix{\chi_C^R:R\ar[r]&\mathrm{Hom}_{R}(C,C)}$ is an isomorphism.

(3)~Ext$^{\geq 1}_{R}(C,C)$=0.

\vskip.25cm From now on, $C$ is a fixed semidualizing $R$-module.

\vskip.25cm\textbf{~1.3} (see \cite{w,hw}) The classes of $C$-projective and $C$-flat modules are defined as
\begin{center}$\mathcal P_{C}(R)=\{C\otimes P|P $ is a projective module\},\\
 $\mathcal F_{C}(R)=\{C\otimes F|F $ is a flat module\}.\end{center}
When $C=R$, we omit the subscript and recover the classes
of projective and flat $R$-modules.

\begin{Lem}\label{L1.1}(\cite[Proposition 5.2]{hw}) Let $\xymatrix{0\ar[r]&W' \ar[r]& W\ar[r]&W''\ar[r]& 0}$ be a short exact sequence of $R$-modules.
 If $W', W''\in \mathcal P_{C}(R)$ (resp. $\mathcal F_{C}(R))$, then $W \in \mathcal P_{C}(R)$ (resp. $\mathcal F_{C}(R))$.
\end{Lem}

\vskip.25cm\textbf{~1.4} (\cite{g04}) Let $\mathcal{X}$ be a class of $R$-modules. A complex $X$ is called an $\mathcal{X}$-complex if $X$ is exact and $Z_i(X)\in \mathcal{X}$ for all $i\in \mathbb{Z}$. We let $\widetilde{\mathcal{X}}$ denote the
class of $\mathcal{X}$-complexes.  A complex $X$ is called projective (resp. flat, $C$-pojective, $C$-flat), if $X$ is a $\mathcal{P}(R)$ (resp. $\mathcal{F}(R)$, $\mathcal{P}_{C}(R)$, $\mathcal{F}_{C}(R)$)-complex.

\vskip.25cm\textbf{~1.5} Let $\mathcal{A}$ be an Abelian category and $\mathcal{B}$ a full subcategory of $\mathcal{A}$. Recall that an exact sequence $\textbf{L}$ in $\mathcal{A}$ is ${\rm Hom}(-,\mathcal{B})$-exact if the sequence ${\rm Hom}(\textbf{L},B)$ is exact for any $B\in \mathcal{B}$.

\vskip.25cm\textbf{~1.6} (\cite[Definition 1.1]{zwl}) An $R$-module $M$ is called $D_{C}$-projective if there exists a $\text{Hom}(-,\mathcal{F}_{C}({R}))$-exact exact sequence \vspace*{-2mm}$$\xymatrix{\cdots\ar[r]^{} & X^{1} \ar[r]^{\delta^{1}} & X^{0}\ar[r]^{\delta^{0}} & X^{-1}\ar[r]^{\delta^{-1}}& X^{- 2}\ar[r]^{\delta^{-2}} &\cdots,}\vspace*{-2mm}$$
of modules with $X^{i}\in \mathcal{P}(R)$ for all $i\geq 0$ and $X^{i}\in \mathcal{P}_{C}(R)$ for all $i<0$ such that $M\cong {\rm Im}\delta^{0}.$

The class of $D_{C}$-projective $R$-modules denoted by $\mathcal{D}_C\mathcal{P}(R)$.

\section{$D_{C}$-projective complexes}\label{sec1}

In this Section, we introduce and study $D_{C}$-projective complexes.

\begin{Def}\label{D2.1} A complex $X$ is called Ding projective with respective to a semidulizing module $C$, simply $D_{C}$-projective, if there exists a ${\rm Hom}(-,\widetilde{\mathcal{F}_C(R)})$-exact exact sequence of complexes
\vspace*{-2mm}$$\xymatrix{\cdots\ar[r]^{} & P_{1} \ar[r]^{f_{1}} & P_{0}\ar[r]^{f_{0}} & Q_{-1}\ar[r]^{f_{-1}}& Q_{- 2}\ar[r]^{f_{-2}} &\cdots}\vspace*{-2mm}$$
with all $P_{i}\in \widetilde{\mathcal{P}(R)}$ and all $Q_{i}\in \widetilde{\mathcal{P}_{C}(R)}$ such that $X\cong {\rm Im}f_{0}$.\end{Def}

By the definition of $D_C$-projective complexes, we have

\begin{Lem}\label{L2.2} If $X$ is a $D_{C}$-projective complex, then ${\rm Ext}^{\geq 1}(X,H)$=0 for any $H\in\widetilde{\mathcal{F}_C(R)}$.\end{Lem}

\begin{Lem}\label{L2.3} Let $X\in{\rm Ch}(R)$. If each $X^{i}\in \mathcal{D}_C\mathcal{P}(R)$ for all $i\in\mathbb{Z}$, then for any $H\in \widetilde{\mathcal{F}_C(R)}$, $\mathcal{H}{\rm om}(X,H)$ if and only if ${\rm Ext}^{1}(X,H)$=0.\end{Lem}
\proof It follows from \cite[Lemma 2.1]{g04}.
\endproof

\begin{Prop}\label{P2.4} Let $X\in{\rm Ch}(R)$. If $X$ is bounded right and each $X^{i}\in \mathcal{D}_C\mathcal{P}(R)$ for all $i\in\mathbb{Z}$, then ${\rm Ext}^{1}(X,H)$=0 for any $H\in \widetilde{\mathcal{F}_C(R)}$.\end{Prop}
\proof
Assume that $\xymatrix{0\ar[r] &H\ar[r]^\mu &G\ar[r]^\nu &X\ar[r] &0}$ is a short exact sequence in ${\rm Ch}(R)$ with $H\in \widetilde{\mathcal{F}_C(R)}$. It suffices to prove that this exact sequence is split.

Without loss of generality, we set
\vspace*{-2mm}$$X=\xymatrix{\cdots\ar[r] &X^2\ar[r]^{\delta^2_X} &X^1\ar[r]^{\delta^1_X} &X^0\ar[r]^{\delta^0_X} &0\ar[r]&0\ar[r]&\cdots.}\vspace*{-2mm}$$
Consider the following commutative diagram
\vspace*{-2mm}$$\xymatrix{        &\vdots \ar[d]^{}                          &\vdots \ar[d]^{}                         &\vdots \ar[d]^{}            \\
 0  \ar[r]^{}       &H^2\ar[d]_{\delta^2_H} \ar[r]^{\mu^2}     &G^2  \ar[d]_{\delta^2_G}  \ar[r]^{\nu^2} & X^2 \ar[d]_{\delta^2_X} \ar[r]^{} &0  \\
 0  \ar[r]^{}       &H^1\ar[d]_{\delta^1_H} \ar[r]^{\mu^1}     &G^1  \ar[d]_{\delta^1_G}  \ar[r]^{\nu^1} & X^1 \ar[d]_{\delta^1_X} \ar[r]^{} &0 \\
 0  \ar[r]^{}       &H^0\ar[d]_{\delta^0_H} \ar[r]^{\mu^0}     &G^0  \ar[d]_{\delta^0_G}  \ar[r]^{\nu^0} & X^0 \ar[d]_{\delta^0_X} \ar[r]^{} &0  \\
 0  \ar[r]^{}       &H^{-1}\ar[d]^{} \ar[r]^{\mu^{-1}}         &G^{-1}  \ar[d]_{}\ar[r]^{\nu^{-1}}       & 0 \ar[d]_{} \ar[r]^{} &0  \\
                    &\vdots                                    &\vdots                                   &  \vdots
                                      }\vspace*{-2mm}$$
Since $H\in \widetilde{\mathcal{F}_C(R)}$, $X^i$ is $D_{C}$-projective, by Lemma \ref{L1.1} and \cite[Proposition 1.4]{zwl},
${\rm Ext}^{1}(X^i,H^i)$=0 for all $i\in \mathbb{Z}$. Thus there exists a $\lambda^i:G^i\longrightarrow H^i$ such that $\lambda^i\mu^i=1_{H^i}$ for all $i\in\mathbb{Z}$. Obviously, $\lambda^i=(\mu^i)^{-1}$  is an isomorphism for all $i<0$. Hence $\delta^i_H\lambda^i=\lambda^{i-1}\delta^i_G$ for all $i<0$.

Note that $\delta_H^{-1}(\lambda^{-1}\delta^0_G-\delta^0_H\lambda^0)=0$, so ${\rm Im}(\lambda^{-1}\delta^0_G-\delta^0_H\lambda^0)\subseteq {\rm Ker}\delta_H^{-1}={\rm Im}\delta_H^{0}$. This implies that
$\lambda^{-1}\delta^0_G-\delta^0_H\lambda^0\in{\rm Hom}(G^0, {\rm Im}\delta^0_H)$. On the other hand, from $(\lambda^{-1}\delta^0_G-\delta^0_H\lambda^0)\mu^0=0$ it follows that ${\rm Ker}\nu^0={\rm Im}\mu^0\subseteq{\rm Ker}(\lambda^{-1}\delta^0_G
-\delta^0_H\lambda^0)$. Thus by the Factor Lemma, there exists $\sigma^0\in{\rm Hom}(X^0, {\rm Im}\delta^0_H)$ such that $\lambda^{-1}\delta^0_G-\delta^0_H\lambda^0=\sigma^0\nu^0$. Since $X^0\in\mathcal{D}_C\mathcal{P}(R)$ and ${\rm Im}\delta^1_H\in \mathcal{F}_C(R)$, ${\rm Ext}^1(X^0,{\rm Im}\delta^1_H)=0$ by \cite[Proposition 1.4]{zwl}. So the sequence
\vspace*{-2mm}$$\xymatrix{0\ar[r]&{\rm Hom}(X^0, {\rm Im}\delta^1_H)\ar[r]&{\rm Hom}(X^0, H^0)\ar[r]&{\rm Hom}(X^0, {\rm Im}\delta^0_H)\ar[r]&0}\vspace*{-2mm}$$ is exact. Hence there exists
a $\tau^0\in{\rm Hom}(X^0, H^0)$ such that $\delta^0_H\tau^0=\sigma^0.$ Take $\omega^0=\tau^0\nu^0+\lambda^0$. Then $\omega^0\in{\rm Hom}(G^0,H^0)$ and
\vspace*{-2mm}$$\delta^0_H\omega^0=\delta^0_H(\tau^0\nu^0+\lambda^0)=\sigma^0\nu^0+\delta^0_H\lambda^0=\lambda^{-1}\delta^0_G,\,\,
\omega^0\mu^0=\tau^0\nu^0\mu^0+\lambda^0\mu^0=1_{H^0}.\vspace*{-2mm}$$

Since $\delta_H^{0}(\omega^0\delta^1_G-\delta^1_H\lambda^1)=0$,
${\rm Im}(\omega^0\delta^1_G-\delta^1_H\lambda^1)\subseteq {\rm Ker}\delta_H^{0}={\rm Im}\delta_H^{1}$.
This implies that $\omega^0\delta^1_G-\delta^1_H\lambda^1\in{\rm Hom}(G^1, {\rm Im}\delta^1_H)$.
On the other hand, since $(\omega^0\delta^1_G-\delta^1_H\lambda^1)\mu^1=0$, one gets ${\rm Ker}\nu^1={\rm Im}\mu^1\subseteq{\rm Ker}(\omega^0\delta^1_G
-\delta^1_H\lambda^1)$. By the Factor Lemma, there exists a $\sigma^1\in{\rm Hom}(X^1, {\rm Im}\delta^1_H)$ such that $\omega^0\delta^1_G-\delta^1_H\lambda^1=\sigma^1\nu^1$. Since ${\rm Ext}^1(X^1,{\rm Im}\delta^2_H)=0$, we have the following exact sequence
\vspace*{-2mm}$$\xymatrix{0\ar[r]&{\rm Hom}(X^1, {\rm Im}\delta^2_H)\ar[r]&{\rm Hom}(X^1, H^1)\ar[r]&{\rm Hom}(X^1, {\rm Im}\delta^1_H)\ar[r]&0.}\vspace*{-2mm}$$
Thus there exists $\tau^1\in{\rm Hom}(X^1, H^1)$ such that $\delta^1_H\tau^1=\sigma^1.$ Put $\omega^1=\tau^1\nu^1+\lambda^1$. Then $\omega^1\in{\rm Hom}(G^1,H^1)$ and
\vspace*{-2mm}$$\delta^1_H\omega^1=\delta^1_H(\tau^1\nu^1+\lambda^1)=\sigma^1\nu^1+\delta^1_H\lambda^1=\omega^0\delta^1_G,\,
\omega^1\mu^1=\tau^1\nu^1\mu^1+\lambda^1\mu^1=1_{H^1}.\vspace*{-2mm}$$

Continue this process, we can get $\omega^i\in{\rm Hom}(G^i,H^i)$ with
$\delta^i_H\omega^i=\omega^{i-1}\delta^i_G$ and $\omega^i\mu^i=1_{H^i}$ for $ i=2,3,\ldots.$ Take $\omega^i=\lambda^i$ when $i<0$ʱ. Then $\omega:G\longrightarrow H$ is a morphism with $\omega\mu=1_H$. Therefore the sequence
\vspace*{-2mm}$$\xymatrix{0\ar[r] &H\ar[r]^\mu &G\ar[r]^\nu &X\ar[r] &0}\vspace*{-2mm}$$ is exact, and so Ext$^{1}(X,H)$=0.
\endproof

Let $\lambda$ be an ordinal number, $(X_\alpha)_{\alpha<\lambda}$ a family subcomplexes of a complex
$X$. Recall that the family $(X_\alpha)_{\alpha<\lambda}$ is a continuous chain of subcomplexes \cite[Definition 2.8]{ei} if $X_{\alpha}\subseteq X_{\beta}$ whenever $\alpha\leq\beta<\lambda$ and if $X_{\beta}=\bigcup_{\alpha<\beta}X_{\alpha}$ whenever $\beta<\lambda$
is a limit ordinal.

\begin{Prop}\label{P2.5} Let $X\in{\rm Ch}(R)$. If $X$ is an exact complex and each $Z_i(X)\in \mathcal{D}_C\mathcal{P}(R)$ for all $i\in\mathbb{Z}$, then ${\rm Ext}^{1}(X,H)$=0 for any $H\in \widetilde{\mathcal{F}_C(R)}$.\end{Prop}
\proof Assume that
\vspace*{-2mm}$$X=\xymatrix{\cdots\ar[r] &X^2\ar[r] &X^1\ar[r] &X^0\ar[r] &X^{-1}\ar[r]&X^{-2}\ar[r]&\cdots}\vspace*{-2mm}$$ is an exact complex and each $Z_i(X)\in \mathcal{D}_C\mathcal{P}(R)$ for all $i\in\mathbb{Z}$. Then by \cite[Theorem 1.12]{zwl}, each $X^i\in \mathcal{D}_C\mathcal{P}(R)$ for all $i\in\mathbb{Z}$.
For any $m\geq0$, let
 \vspace*{-2mm}$$X_m=\xymatrix{\cdots\ar[r] &X^0\ar[r] &X^{-1}\ar[r]&\cdots\ar[r]&X^{-m}\ar[r]&Z_{-m-1}(X)\ar[r]&0.}\vspace*{-2mm}$$
Then $(X_m)_{m\geq0}$ is a continuous chain of subcomplexes of $X$ and $X=\bigcup_{m\geq0}X_{m}$. Since $X_{m+1}/X_m\cong \overline{Z_{-m-1}(X)}[-m-2]$ for $m=0,1,\ldots,$ by Proposition \ref{P2.4}, Ext$^{1}(X_0,H)=0$,
Ext$^{1}(X_{m+1}/X_m,H)=0$ for any $ m=0,1,\ldots$ and any $H\in \widetilde{\mathcal{F}_C(R)}$. Therefore Ext$^{1}(X,H)=0$ for any $H\in \widetilde{\mathcal{F}_C(R)}$ by \cite[Theorem 2.9]{ei}.\endproof

\begin{Cor} \label{C2.6} If $X\in \widetilde{\mathcal{P}_C(R)}$, then $\mathcal{H}{\rm om}(X,H)$ is exact for any $H\in \widetilde{\mathcal{F}_C(R)}$.\end{Cor}
\proof Since $X\in \widetilde{\mathcal{P}_C(R)}$, $X$ is an exact complex, and each $Z_i(X)\in \mathcal{D}_C\mathcal{P}(R)$ for any $i\in\mathbb{Z}$ by \cite[Proposition 1.8]{zwl}. Thus Ext$^{1}(X,H)$=0 for any $H\in \widetilde{\mathcal{F}_C(R)}$ by Proposition \ref{P2.5}. On the other hand, each $X^i\in \mathcal{D}_C\mathcal{P}(R)$ by \cite[Theorem 1.12]{zwl}. Thus $\text{$\mathcal{H}$om}(X,H)$ is exact for any $H\in \widetilde{\mathcal{F}_C(R)}$ by Lemma \ref{L2.3}.\endproof

\begin{Lem}\label{L2.7}
Let $\xymatrix{0\ar[r] &M\ar[r]^f&L\ar[r]&N\ar[r]&0}$ be an exact sequence of $R$-modules. If $N\in \mathcal{D}_C\mathcal{P}(R)$, $L\in \mathcal{P}_C(R)$, then for any $L'\in\mathcal{P}_C(R)$ and any $f'\in{\rm Hom}(M, L')$, ${\rm Coker}\alpha\in \mathcal{D}_C\mathcal{P}(R)$ where $\xymatrix{\alpha=(f,f'):M\ar[r]&L\oplus L'.}$
\end{Lem}
\proof
Suppose that $L'\in\mathcal{P}_C(R)$, $f'\in\text{Hom}(M, L')$ and $\xymatrix{\alpha=(f,f'):M\ar[r]&L\oplus L'.}$
Then we have the following exact sequence
 \vspace*{-2mm}$$\xymatrix{0\ar[r]&M\ar[r]^{\alpha}  &L\oplus L'\ar[r]&{\rm Coker}\alpha\ar[r]&0.}\vspace*{-2mm}$$
By the Factor Lemma, there exists an epimorphism $\mu:{\rm Coker}\xymatrix{\alpha\ar[r]& N}$ such that the following diagram commutes
  \vspace*{-2mm}$$\xymatrix{
   0\ar[r]_{ }       &M\ar@{=}[d] \ar[r]^{\alpha}     &L\oplus L'  \ar[d]_{\pi}  \ar[r]_{ } & \textrm{Coker}\alpha \ar[d]_{\mu} \ar[r]_{ } &0  \\0 \ar[r]_{ }     &M\ar[r]^{f}  &L\ar[r]_{ }   &N\ar[r]_{ }  &0 }\vspace*{-2mm}$$
where $\pi:\xymatrix{L\oplus L'\ar[r] &L}$ is the canonical projection. By the Snake Lemma, ${\rm Ker}\mu\cong L'$. Thus ${\rm Ker}\mu\in\mathcal{P}_C(R)$. So ${\rm Ker}\mu\in \mathcal{D}_C\mathcal{P}(R)$ by \cite[Proposition 1.8]{zwl}. Hence ${\rm Coker}\alpha\in \mathcal{D}_C\mathcal{P}(R)$ by \cite[Proposition 1.10]{zwl}.
\endproof

Now, we can achieve a characterization of $D_{C}$-projective complexes.

\begin{The}\label{T2.8}  Let $X\in{\rm Ch}(R)$. Then $X$ is a $D_{C}$-projective complex if and only if each $X^i\in \mathcal{D}_C\mathcal{P}(R)$ for all $i\in\mathbb{Z}$ and $\mathcal{H}{\rm om}(X,H)$ is exact for any $H\in \widetilde{\mathcal{F}_C(R)}$.\end{The}

\proof $\Longrightarrow)$ Assume that ${X=\xymatrix{\cdots\ar[r] &X^{i+1}\ar[r]& X^{i}\ar[r]& X^{i-1}\ar[r]& \cdots}}$ is a $D_{C}$-projective complex. Then there exists a ${\rm Hom}(-,\widetilde{\mathcal{F}_C(R)})$-exact exact sequence of complexes
\vspace*{-2mm}$$\xymatrix{\mathbf{E}=\cdots\ar[r]^{} & P_{2} \ar[r]^{f_{2}}& P_{1} \ar[r]^{f_{1}} & P_{0}\ar[r]^{f_{0}} & Q_{-1}\ar[r]^{f_{-1}}& Q_{- 2}\ar[r]^{f_{-2}} &\cdots}\vspace*{-2mm}$$ with each $P_{i}\in \widetilde{\mathcal{P}(R)}$, each $Q_{i}\in\widetilde{\mathcal{P}_C(R)}$ such that
$X\cong{\rm Im}f_0$. Of course for any $i\in \mathbb{Z}$, we have the following exact sequence of $R$-modules
\vspace*{-2mm}$$\xymatrix{\mathbf{E}^i=\cdots\ar[r]^{} & P_{2}^{i} \ar[r]^{f_{2}^i}& P_{1}^{i} \ar[r]^{f_{1}^i} & P_{0}^{i}\ar[r]^{f_{0}^i} & Q_{-1}^{i}\ar[r]^{f_{-1}^i}& Q_{- 2}^{i}\ar[r]^{f_{-2}^i} &\cdots,}\vspace*{-2mm}$$ and it does have $X^{i}\cong{\rm Im}f_0^i$, where all $P_{j}^{i}\in \mathcal{P}(R)$ and all $Q_{j}^{i}\in \mathcal{P}_C(R)$ (by Lemma \ref{L1.1}). Let $F\in\mathcal{F}_C(R)$. Then we have the following commutative diagram with the the top row exact
\vspace*{-2mm}$$\xymatrix{
\cdots\ar[r]&\text{Hom}(Q_{-1}, \overline {F}[i]) \ar[d]_{} \ar[r]_{ }     &\text{Hom}(P_{0}, \overline{F} [i])  \ar[d]_{}\ar[r]_{ } & \text{Hom}(P_{1}, \overline{F }[i])  \ar[d]_{}\ar[r]&\cdots\\\cdots \ar[r]&\text{Hom}({Q_{-1}^i}, {F}) \ar[r]_{ }     &\text{Hom}({P_{0}^i},{F})  \ar[r]_{ } & \text{Hom}({P_{1}^i}, {F})\ar[r]&\cdots }\vspace*{-2mm}$$
Now the bottom row is exact since the vertical maps are all isomorphism by \cite[Lemma 3.1]{g04}. This shows that $\mathbf{E}^i$ remains exact after applying $\text{Hom}(-,F)$ for any $F\in \mathcal{F}_C(R)$. Thus each $X^{i}\in\mathcal{D}_C\mathcal{P}(R)$ for all $i\in \mathbb{Z}$. The remainder follows from Lemma \ref{L2.2} and Lemma \ref{L2.3}.

$\Longleftarrow)$ Since $X^{i}\in \mathcal{D}_{C}\mathcal{P}(R)$, there exists an exact sequence of $R$-modules
\vspace*{-2mm}$$\xymatrix{ 0\ar[r]^{} & X^{i} \ar[r]^{f^{i}} & Q^{i}\ar[r]^{} & Y^{i}\ar[r]^{}&0}\vspace*{-2mm}$$ such that $Q^{i}\in \mathcal{P}_{C}(R)$ and $Y^{i}\in \mathcal{D}_{C}\mathcal{P}(R)$ by \cite[Proposition 1.13]{zwl} for all $i\in \mathbb{Z}$.
For any $i\in \mathbb{Z}$ and any $(x,y)\in Q^{i}\oplus Q^{i-1}$, define $\delta^{i}(x,y)=(y,0)$. Then
\vspace*{-2mm}$$\xymatrix{Q_{-1}=\cdots \ar[r]_{ } &Q^{i}\oplus Q^{i-1}\ar[r]^{\delta^{i} } &Q^{i-1}\oplus Q^{i-2}\ar[r]^{\delta^{i-1} } &Q^{i-2}\oplus Q^{i-3}\ar[r]^{ }&\cdots }\vspace*{-2mm}$$
is a $C$-projective complex. Let $i\in \mathbb{Z}$, put $\alpha_{i}=(f^{i},f^{i-1}\sigma^{i}):\xymatrix{X^{i}\ar[r]& Q^{i}\oplus Q^{i-1}.}$ Then $\alpha=(\alpha_i):\xymatrix{X\ar[r]&Q_{-1}}$ is a morphism between the following two complexes
\vspace*{-2mm}$$\xymatrix{
X=\cdots\ar[r]_{ } &X^{i+1}\ar[d]_{\alpha_{i+1}}  \ar[r]^{\sigma^{i+1}} &X^{i}\ar[d]_{\alpha_{i}}  \ar[r]^{\sigma^{i}}&X^{i-1}\ar[d]_{\alpha_{i-1}}\ar[r]^{ } &\cdots \\Q_{-1}=\cdots \ar[r]_{ } &Q^{i+1}\oplus Q^{i}\ar[r]^{\delta^{i+1} } &Q^{i}\oplus Q^{i-1}\ar[r]^{\delta^{i} } &Q^{i-1}\oplus Q^{i-2}\ar[r]^{ }&\cdots }\vspace*{-2mm}$$
So we get an exact sequence of complexes
\vspace*{-2mm}$$\xymatrix{ 0\ar[r]^{} & X\ar[r]^{\alpha} & Q_{-1}\ar[r]^{} & K_{-1}\ar[r]^{}&0,}\vspace*{-2mm}$$
where $K_{-1}^{i}={\rm Coker}\alpha_{i}, \forall i\in \mathbb{Z}$. It follows from Lemma \ref{L2.7} that each $K_{-1}^{i}\in \mathcal{D}_{C}\mathcal{P}(R)$. Thus the
sequence of complexes
\vspace*{-2mm}$$\xymatrix{ 0\ar[r]_{ } &\text{$\mathcal{H}$om}(K_{-1},H) \ar[r]_{ }     &\text{$\mathcal{H}$om}(Q_{-1},H)  \ar[r]_{ } &\text{$\mathcal{H}$om}(X,H)\ar[r]_{ }&0 }\vspace*{-2mm}$$ is exact for any $H\in \widetilde{\mathcal{F}_C(R)}$ by Lemma \cite[Proposiiton 1.4]{zwl}.
Since $\text{$\mathcal{H}$om}(X,H)$ is exact and $\text{$\mathcal{H}$om}(Q_{-1},H)$ is exact by Corollary \ref{C2.6}, $\text{$\mathcal{H}$om}(K_{-1},H)$ is exact. Hence Ext$^{1}(K_{-1},H)$=0 by Lemma \ref{L2.3}. This yields the following exact sequence
\vspace*{-2mm}$$\xymatrix{ 0\ar[r]_{ }       &\text{Hom}(K_{-1},H)\ar[r]_{ }     &\text{Hom}(Q_{-1},H)  \ar[r]_{ } &\text{Hom}(X,H )\ar[r]_{ }&0 .}\vspace*{-2mm}$$
Note that $K_{-1}$ has the same properties
as $X$, we may use the same procedure to construct a ${\rm Hom}(-,\widetilde{\mathcal{F}_C(R)})$-exact exact sequence of complexes
\vspace*{-2mm}$$\xymatrix{0\ar[r]& X\ar[r]& Q_{-1}\ar[r]& Q_{-2} \ar[r]&\cdots} \eqno(\dagger)\vspace*{-2mm}$$
with each $Q_{i}\in \widetilde{\mathcal{P}_C(R)}$ for any $i\in \mathbb{Z}$.

Take a projective resolution of $X$
\vspace*{-2mm}$$\xymatrix{\cdots\ar[r]& P_{1}\ar[r]^{f_1}&P_{0}\ar[r]^{f_0}&X\ar[r]& 0.} \eqno(\ddagger)\vspace*{-2mm}$$
Set $K_j=\text{Ker}f_j,\, j=0,1,2,\ldots.$ Then $K_j^i\in \mathcal{D}_{C}\mathcal{P}(R)$ by \cite[Theorem 1.12]{zwl} for any $j\geq0$ and any $i\in \mathbb{Z}$. By the hypothesis and Lemma \ref{L2.3}, the sequence
\vspace*{-2mm}$$\xymatrix{ 0\ar[r]_{ }       &\text{Hom}(X,H)\ar[r]_{ }     &\text{Hom}(P_0,H)  \ar[r]_{ } &\text{Hom}(K_0,H )\ar[r]_{ }&0 }\vspace*{-2mm}$$
is exact for any $H\in \widetilde{\mathcal{F}_C(R)}$. On the other hand, by Lemma \cite[Proposition 1.4]{zwl}, the sequence of complexes
\vspace*{-2mm}$$\xymatrix{ 0\ar[r]_{ }       &\mathcal{H}\text{om}(X,H)\ar[r]_{ }     &\mathcal{H}\text{om}(P_0,H)  \ar[r]_{ } &\mathcal{H}\text{om}(K_0,H )\ar[r]_{ }&0 }\vspace*{-2mm}$$ is exact for any $H\in \widetilde{\mathcal{F}_C(R)}$ since each $X^i\in \mathcal{D}_{C}\mathcal{P}(R)$. Then $\mathcal{H}\text{om}(K_0,H )$ is exact since both $\mathcal{H}\text{om}(X,H)$ and $\mathcal{H}\text{om}(P_0,H)$ are all exact. Continue this process one can prove that $(\ddagger)$ is ${\rm Hom}(-,\widetilde{\mathcal{F}_C(R)})$-exact.

Now, assembling the sequence $(\dagger)$ and $(\ddagger)$ together, we get a ${\rm Hom}(-,\widetilde{\mathcal{F}_C(R)})$-exact exact sequence of complexes
\vspace*{-2mm}$$\xymatrix{\mathbf{E}=\cdots\ar[r]^{} & P_{1} \ar[r]^{} & P_{0}\ar[r]^{} & Q_{-1}\ar[r]^{}& Q_{- 2}\ar[r]^{} &\cdots}\vspace*{-2mm}$$ with all $P_{i}\in \widetilde{\mathcal{P}(R)}$ and all $Q_{i}\in \widetilde{\mathcal{P}_C(R)}$ such that $X\cong$ Im$(\xymatrix{P_{0} \ar[r]& Q_{-1}).}$ Therefore, $X$ is a $D_C$-projective complex.\endproof

\begin{Cor}\label{C2.9} Projective complexes and $C$-projective complexes are $D_{C}$-projective complexes.\end{Cor}
\proof It follows from Theorem \ref{T2.8}, \cite[Proposition 1.8]{zwl} and Corollary \ref{C2.6}.\endproof

Let $\mathcal{A}$ be an Abelian category. According to \cite{h}, a class $\mathcal{X}$ of objects of $\mathcal{A}$ is said to be projectively resolving if all projective objects of $\mathcal{A}$ are contained in $\mathcal{X}$ and for every short exact sequence $\xymatrix{0\ar[r]&X'\ar[r]&X\ar[r]&X''\ar[r]&0}$ in $\mathcal{A}$, if $X''\in \mathcal{X}$, then $X\in \mathcal{X}$ if and only if $X'\in \mathcal{X}$.

\begin{Cor}\label{C2.10} The class of $D_{C}$-projective complexes is projectively resolving.\end{Cor}
\proof By Corollary \ref{C2.9}, every projective complex is $D_{C}$-projective. Now consider an exact sequence in ${\rm Ch}(R)$
\vspace*{-2mm} $$\xymatrix{0\ar[r] &X\ar[r] &Y\ar[r] &Z\ar[r]&0}\vspace*{-2mm}$$
with $Z$ $D_{C}$-projective. Then $Z^{i}\in \mathcal{D}_{C}\mathcal{P}(R)$ for any $i\in\mathbb{Z}$ and $\text{$\mathcal{H}$om}(Z, H)$ is exact for any $H\in \widetilde{\mathcal{F}_C(R)}$ by Theorem \ref{T2.8}. So the sequence \vspace*{-2mm}$$\xymatrix{ 0\ar[r]&\text{$\mathcal{H}$om}(Z,H)\ar[r]     &\text{$\mathcal{H}$om}(Y, H) \ar[r] &\text{$\mathcal{H}$om}(X, H)\ar[r]_{ }&0 }\vspace*{-2mm}$$ is exact for any $H\in \widetilde{\mathcal{F}_C(R)}$ by \cite[Proposition 1.4]{zwl}. Thus if $X$ is a $D_{C}$-projective complex, then $X^{i}\in \mathcal{D}_{C}\mathcal{P}(R)$ for any $i\in\mathbb{Z}$ and $\text{$\mathcal{H}$om}(X,H)$ is exact for any $H\in \widetilde{\mathcal{F}_C(R)}$ by Theorem \ref{T2.8}. Hence $\text{$\mathcal{H}$om}(Y,H)$ is exact for any $H\in \widetilde{\mathcal{F}_C(R)}$, and by \cite[Theorem 1.12]{zwl}, $Y^{i}\in\mathcal{D}_{C}\mathcal{P}(R)$ for any $i\in\mathbb{Z}$. Therefore $Y$ is a $D_{C}$-projective complex by Theorem \ref{T2.8}.
The case $Y$ is a $D_{C}$-projective complex can be proved similarly.\endproof

\begin{Cor}\label{C2.11} The class of $D_{C}$-projective complexes is closed under
direct summands and direct sums.\end{Cor}
\proof Suppose that $X$ is a $D_{C}$-projective complex and $A\oplus B=X$. Then $A^{i}\in \mathcal{D}_{C}\mathcal{P}(R)$ for any
$i\in \mathbb{Z}$ by Theorem \ref{T2.8} and \cite[Theorem 1.12]{zwl}. Let $H\in \widetilde{\mathcal{F}_C(R)}$. Then $\text{$\mathcal{H}$om}(X,H)$ is exact by Theorem \ref{T2.8}. Thus $\text{$\mathcal{H}$om}(A, H)$ is exact since $\text{$\mathcal{H}$om}(A,H)\oplus \text{$\mathcal{H}$om}(B, H)\cong \text{$\mathcal{H}$om}(A\oplus B,H)$. Hence $A$ is a $D_{C}$-projective complex by Theorem \ref{T2.8}.

Let $\{X_\lambda\}_{\lambda\in \Lambda}$ be a collection of $D_{C}$-projective complexes. Then $\bigoplus_{\lambda\in\Lambda}X_\lambda^{i} \in \mathcal{D}_{C}\mathcal{P}(R)$ for any
$i\in \mathbb{Z}$ by Theorem \ref{T2.8} and \cite[Proposition 1.11]{zwl}. Let $H\in \widetilde{\mathcal{F}_C(R)}$. Then $\text{$\mathcal{H}$om}(X_{\lambda},H)$ is exact for any $\lambda\in\Lambda$. Since $\text{$\mathcal{H}$om}(\bigoplus_{\lambda\in \Lambda}X_{\lambda}, H)\cong \prod_{\lambda\in\Lambda}\text{$\mathcal{H}$om}(X_{\lambda},H)$, $\text{$\mathcal{H}$om}(\bigoplus_{\lambda\in \Lambda}X_{\lambda}, H)$
is exact. So $\bigoplus_{\lambda\in\Lambda}X_{\lambda}$ is a $D_{C}$-projective complex by Theorem \ref{T2.8}.\endproof

\begin{Cor}\label{C2.12} Let $\xymatrix{0\ar[r] &X\ar[r] &Y\ar[r] &Z\ar[r]&0}$ be a short exact sequence of complexes.
If $X,Y$ are $D_{C}$-projective, then the following conditions are equivalent:

  (1) $Z$ is a $D_{C}$-projective complex.

  (2) $Z^{i}\in \mathcal{D}_{C}\mathcal{P}(R)$ for all $ i\in \mathbb{Z}$.

  (3) ${\rm Ext}^{1}(Z,H)$=0 for any $H\in \widetilde{\mathcal{F}_C(R)}$.\end{Cor}

\proof (1)$\Longrightarrow$(3) It follows from Lemma \ref{L2.2}.

(3)$\Longrightarrow$(2) Let $i\in \mathbb{Z}$. Consider the exact sequence of $R$-modules
\vspace*{-2mm} $$\xymatrix{0\ar[r]& X^{i}\ar[r]&Y^{i}\ar[r]&Z^{i}\ar[r]&0 .}\vspace*{-2mm}$$
By Theorem \ref{T2.8}, $X^{i},Y^{i}\in \mathcal{D}_{C}\mathcal{P}(R)$. It suffices to show that Ext$^{1}(Z^{i},F)=0$ for any $F\in \mathcal{F}_C(R)$ by \cite[Corollary 1.15]{zwl}.

Let $F\in \mathcal{F}_C(R)$. Then $\overline{F}[i]\in \widetilde{\mathcal{F}_C(R)}$. Thus Ext$^{1}(Z,\overline{F}[i])$=0 by (3).
Hence Ext$^{1}(Z^{i},F)$=0 since Ext$^{1}(Z^{i},F)\cong$ Ext$^{1}(Z,\overline{F}[i])$ by \cite[Lemma 3.1]{g04}.

(2)$\Longrightarrow$(1) Assume that $H\in \widetilde{\mathcal{F}_C(R)}$. Since each $Z^{i}\in \mathcal{D}_{C}\mathcal{P}(R)$, the sequence
\vspace*{-2mm} $$\xymatrix{ 0\ar[r]_{ }       &\text{$\mathcal{H}$om}(Z, H)\ar[r]_{ }     &\text{$\mathcal{H}$om}(Y, H) \ar[r]_{ } &\text{$\mathcal{H}$om}(X, H)\ar[r]_{ }&0 }\vspace*{-2mm}$$ is exact by \cite[Proposition 1.4]{zwl}. Then $\text{$\mathcal{H}$om}(Z,H)$ is exact since $\text{$\mathcal{H}$om}(X, H)$ and $\text{$\mathcal{H}$om}(Y, H)$ are all exact. So $Z$ is a $D_{C}$-projective complex by Theorem \ref{T2.8}.\endproof

\begin{Cor}\label{C2.13} For any ${\rm Hom}(-,\widetilde{\mathcal{F}_C(R)})$-exact exact sequence\vspace*{-2mm}$$\xymatrix{
\mathbf{E}=\cdots\ar[r]^{} & P_{1} \ar[r]^{f_{1}} & P_{0}\ar[r]^{f_{0}} & Q_{-1}\ar[r]^{f_{-1}}& Q_{- 2}\ar[r]^{f_{-2}} &\cdots}\vspace*{-2mm}$$  Æäwith all $P_{i}\in \widetilde{\mathcal{P}(R)}$ and all $Q_{i}\in \widetilde{\mathcal{P}_C(R)}$, each ${\rm Coker}f_{i}$ is a $D_{C}$-projective complex for any $i\in \mathbb{Z}$.\end{Cor}

\proof Let $i\in \mathbb{Z}$, set $N_{i}={\rm Coker}f_{i+1}$. Then $N_{0}$ is $D_{C}$-projective by the definition of $D_{C}$-projective complexes. Thus $N_i$ is $D_C$-projective for any $i\geq0$ by Corollaries \ref{C2.9} and \ref{C2.10}. So it left to show that $N_i$ is $D_C$-projective for any $i<0$.

For any $m\in \mathbb{Z}$, we have a ${\rm Hom}(-,\mathcal{F}_C(R))$-exact exact sequence of $R$-modules\vspace*{-2mm}$$\xymatrix{
\mathbf{E}^m=\cdots\ar[r]^{} & P_{1}^{m} \ar[r]^{f_{1}^{m}} & P_{0}^{m}\ar[r]^{f_{0}^{m}} & Q_{-1}^{m}\ar[r]^{f_{-1}^{m}}& Q_{- 2}^{m}\ar[r]^{f_{-2}^{m}} &\cdots,}\vspace*{-2mm}$$ where all $P_i^m\in \mathcal{P}(R)$ and all $Q_{i}^m\in \mathcal{P}_C(R)$. Then each $N_{i}^m={\rm Coker}f_{i+1}^m\in \mathcal{D}_C\mathcal{P}(R)$ for any $i\in \mathbb{Z}$ by \cite[Proposition 1.13]{zwl}. 
Thus $N_i$ is $D_C$-projective for any $i<0$ inductively by Corollary \ref{C2.12}.\endproof

The next two Lemmas play a crucial role in the rest of our discussion.

\begin{Lem}\label{L2.14} Let \vspace*{-2mm}
$$\xymatrix{0\ar[r]_{}&A\ar[r]_{}& G_{1}\ar[r]^{f}&G_{0}\ar[r]^{}& X\ar[r]_{}&0}\eqno(2.1)\vspace*{-2mm}$$
be an exact sequence in ${\rm Ch}(R)$ with $G_{0}, G_{1}$ $D_{C}$-projective. Then

(1) We have the following exact sequences\vspace*{-2mm}$$\xymatrix{ 0\ar[r]&A\ar[r]&Q\ar[r]&G\ar[r]&X\ar[r]&0}\eqno(2.2)\vspace*{-2mm}$$
and \vspace*{-2mm}$$\xymatrix{ 0\ar[r]& A\ar[r]& W\ar[r]& P\ar[r]&X\ar[r]&0}\eqno(2.3)\vspace*{-2mm}$$
with $Q\in \widetilde{\mathcal{P}_C(R)}$, $P\in \widetilde{\mathcal{P}(R)}$ and $G, W$ $D_{C}$-projective.

(2) If the sequence (2.1) is ${\rm Hom}(-,\widetilde{\mathcal{F}_C(R)})$-exact,
then so are (2.2) and (2.3).
\end{Lem}
\proof (1) Since $G_{1}$ is $D_{C}$-projective, there exists an exact sequence
\vspace*{-2mm}$$\xymatrix{0\ar[r]&G_{1}\ar[r]&Q\ar[r]&G_{-1}\ar[r]&0}\vspace*{-2mm}$$
with $Q\in \widetilde{\mathcal{P}_C(R)}$ and $G_{-1}$ $D_{C}$-projective by the definition of $D_C$-projective complexes and Corollary \ref{C2.13}. Then we have the following pushout diagram\vspace*{-2mm}$$\xymatrix{
                          &                                        & 0 \ar[d]_{}                  & 0  \ar[d]_{}             &  \\
 0  \ar[r]^{}       &A     \ar@{=}[d]^{} \ar[r]^{}     &G_{1}  \ar[d]_{}  \ar[r]^{ } & {\rm Im}f \ar[d]_{} \ar[r]^{ } &0  \\
 0 \ar[r]^{}        &A\ar[r]^{}  &Q \ar[d]^{ } \ar[r]^{}   &B \ar[d]^{ } \ar[r]^{}  &0  \\
                                       &                 & G_{-1} \ar[d]^{} \ar@{=}[r]        &G_{-1} \ar[d]^{} &\\
                                       &                 & 0                            &0   }\vspace*{-2mm}$$
Consider the following pushout diagram
\vspace*{-2mm}$$\xymatrix{
                                                                  & 0 \ar[d]_{}                  & 0  \ar[d]_{}             &  &  \\
 0  \ar[r]^{}       &{\rm Im}f     \ar[d]^{} \ar[r]^{}     &G_{0}  \ar[d]_{}  \ar[r]^{ } & X \ar@{=}[d] \ar[r]^{ } &0  \\
 0 \ar[r]^{}        &B\ar[d]^{ }\ar[r]^{}  &G\ar[d]^{ } \ar[r]^{}   &X\ar[r]^{}  &0  \\
                                                      & G_{-1} \ar[d]^{} \ar@{=}[r]        &G_{-1}  \ar[d]^{} &  &\\
                                                       & 0                            &0.   }\vspace*{-2mm}$$
Connecting the middle rows in the above two diagrams, we get the exact sequence\vspace*{-2mm}$$\xymatrix{0\ar[r]&A\ar[r]& Q\ar[r]&G\ar[r]& X\ar[r]&0.}\vspace*{-2mm}$$
Since $G_{0},G_{-1}$ are $D_{C}$-projective, then so is $G$ by Corollary \ref{C2.10}. Now the first desired exact sequence (2.2) follows.

Since $G_{0}$ is $D_{C}$-projective, there exists an exact sequence\vspace*{-2mm}$$\xymatrix{0\ar[r]&G^{1}\ar[r]&P\ar[r]&G_{0}\ar[r]&0}\vspace*{-2mm}$$
with $P\in \widetilde{\mathcal{P}(R)}$ and $G^{1}$ $D_{C}$-projective. Then we have the following pullback diagram\vspace*{-2mm}$$\xymatrix{
                                                                  & 0 \ar[d]_{}                  & 0  \ar[d]_{}             &                                     &  \\
      &G^{1}    \ar[d]^{} \ar@{=}[r]    &G^{1}   \ar[d]_{}   &   \\
 0 \ar[r]^{}        &N\ar[d]^{}\ar[r]^{}  &P \ar[d]^{ } \ar[r]^{}   &X\ar@{=}[d] \ar[r]^{}  &0  \\
                                    0 \ar[r]^{}                &{\rm Im}f \ar[d]^{} \ar[r]^{}       &G_{0} \ar[d]^{} \ar[r]^{}&X\ar[r]^{}&0\\
                                                      & 0                            &0  }\vspace*{-2mm}$$
And consider the following pullback diagram
\vspace*{-2mm}$$\xymatrix{
                                                       &           & 0 \ar[d]_{}                  & 0  \ar[d]_{}             &                                      \\
 &     &G^{1}    \ar[d]^{} \ar@{=}[r]    &G^{1}   \ar[d]_{}   &   &\\
 0 \ar[r]^{}        &A\ar@{=}[d]\ar[r]^{}  &W \ar[d]^{ } \ar[r]^{}   &N\ar[d]^{} \ar[r]^{}  &0  \\
                                    0 \ar[r]^{}                &A  \ar[r]^{}       &G_{1} \ar[d]^{} \ar[r]^{}&{\rm Im}f\ar[d]^{} \ar[r]^{}&0\\
                                               &       & 0                            &0   }\vspace*{-2mm}$$
Since both $G^1,G_1$ are $D_C$-projective, so is $W$ by Corollary \ref{C2.10}. Connecting
the middle rows in the above two diagrams, then we get the second desired exact
sequence (2.3).

(2) Let $H\in \widetilde{\mathcal{F}_C(R)}$. Note that ${\rm Ext}^{i\geq1}
(X,H)=0$ for any $D_C$-projective complex $X$ by Lemma \ref{L2.2}. If the exact sequence (2.1) is ${\rm Hom}(-, \widetilde{\mathcal{F}_C(R)})$-exact, then Ext$^{1}({\rm Im}f,H)$=0 and Ext$^{1}(X,H)={\rm Ext}^{2}(X,H)=0.$ So in the proof of
(1), Ext$^{1}(B,H)$=0 and Ext$^{1}(N,H)$=0. Thus the exact sequences (2.2) and
(2.3) are ${\rm Hom}(-, \widetilde{\mathcal{F}_C(R)})$-exact.\endproof

\begin{Lem}\label{L2.15} Let $n\geq1$ and  \vspace*{-2mm}
$$\xymatrix{0\ar[r]_{}&A\ar[r]& G_{n-1}\ar[r]&\cdots\ar[r]&G_{1}\ar[r]&G_{0}\ar[r]& X\ar[r]&0}\eqno(2.4)\vspace*{-2mm}$$
an exact sequence in ${\rm Ch}(R)$ with all $G_{i}$ $D_{C}$-projective. Then

(1) There exist exact sequences
\vspace*{-2mm}$$\xymatrix{0\ar[r]_{}&A\ar[r]& Q_{n-1}\ar[r]&\cdots\ar[r]&Q_{1}\ar[r]&Q_{0}\ar[r]& Y\ar[r]&0}\eqno(2.5)\vspace*{-2mm}$$
and
\vspace*{-2mm}$$\xymatrix{ 0\ar[r]& X\ar[r]& Y\ar[r]& U\ar[r]&0}\vspace*{-2mm}$$
in ${\rm Ch}(R)$ with all $Q_i\in \widetilde{\mathcal{P}_C(R)}$ and $U$ $D_C$-projective.

(2) There exist exact sequences
\vspace*{-2mm}$$\xymatrix{0\ar[r]&B\ar[r]& P_{n-1}\ar[r]&\cdots\ar[r]&P_{1}\ar[r]&P_{0}\ar[r]& X\ar[r]&0}\eqno(2.6)\vspace*{-2mm}$$
and
\vspace*{-2mm}$$\xymatrix{ 0\ar[r]& V\ar[r]& B\ar[r]& A\ar[r]&0}\vspace*{-2mm}$$
in ${\rm Ch}(R)$ with all $P_i\in \widetilde{\mathcal{P}(R)}$ and $V$ $D_C$-projective.

(3) If the exact sequence (2.4) is ${\rm Hom}(-,\widetilde{\mathcal{F}_C(R)})$-exact, then so are (2.5) and (2.6).\end{Lem}

\proof  (1) We proceed by induction on $n$.

When $n=1$, we have an exact sequence $\xymatrix{0\ar[r]&A\ar[r]&G_0\ar[r]&X\ar[r]&0}$ in ${\rm Ch}(R)$. Since $G_0$ is $D_C$-projective, we have a ${\rm Hom}(-, \widetilde{\mathcal{F}_C(R)})$-exact exact
sequence $\xymatrix{0\ar[r]&G_0\ar[r]&Q_0\ar[r]&U\ar[r]&0}$ with $Q_0\in\widetilde{\mathcal{P}_C(R)}$ and $U\,D_C$-projective. Consider the following pushout diagram
\vspace*{-2mm}$$\xymatrix{
                          &                                        & 0 \ar[d]_{}                  & 0  \ar[d]_{}             &  \\
 0  \ar[r]^{}       &A     \ar@{=}[d]^{} \ar[r]^{}     &G_{0}  \ar[d]_{}  \ar[r]^{ } &X \ar[d]_{} \ar[r]^{ } &0  \\
 0 \ar[r]^{}        &A\ar[r]^{}  &Q_0 \ar[d]^{ } \ar[r]^{}   &Y \ar[d]^{ } \ar[r]^{}  &0  \\
                                       &                 & U \ar[d]^{} \ar@{=}[r]        &U \ar[d]^{} &\\
                                       &                 & 0                            &0.     }\vspace*{-2mm}$$
The middle row and the last column in the above diagram are the desired two exact sequences.

Now assume that $n\geq 2$ and we have an exact sequence \vspace*{-2mm}
$$\xymatrix{0\ar[r]_{}&A\ar[r]& G_{n-1}\ar[r]&G_{n-2}\ar[r]&\cdots\ar[r]&G_{1}\ar[r]&G_{0}\ar[r]& X\ar[r]&0 }\vspace*{-2mm}$$
in ${\rm Ch}(R)$ with all $G_{i}$ $D_{C}$-projective. Put $\xymatrix{K={\rm Coker}(G_{n-1}\ar[r]&G_{n-2}).}$ By Lemma \ref{L2.14}, we get an exact sequence
\vspace*{-2mm}$$\xymatrix{0\ar[r]_{}&A\ar[r]& Q_{n-1}\ar[r]&G_{n-2}'\ar[r]&K\ar[r]&0}\eqno(2.7)\vspace*{-2mm}$$
in ${\rm Ch}(R)$ with $Q_{n-1}\in \widetilde{\mathcal{P}_C(R)}$ and $G_{n-2}'$ $D_{C}$-projective. Set $\xymatrix@-0.5pc{A'={\rm Im}(Q_{n-1}\ar[r]&G_{n-2}').}$
Then we get an exact sequence
\vspace*{-2mm}
$$\xymatrix{0\ar[r]_{}&A'\ar[r]& G_{n-2}'\ar[r]& G_{n-3}\ar[r]&\cdots\ar[r]&G_{1}\ar[r]&G_{0}\ar[r]& X\ar[r]&0}\vspace*{-2mm}$$in ${\rm Ch}(R)$.
Now we get the assertion by the induction hypothesis.

(2) The proof is dual to that of (1).

(3) If the exact sequence (2.4) is $\text{Hom}(-,\widetilde{\mathcal{F}_C(R)})$-exact, then the middle rows in the
above commutative diagram is also $\text{Hom}(-,\widetilde{\mathcal{F}_C(R)})$-exact. On the other hand, we
can choose (2.7) to be $\text{Hom}(-,\widetilde{\mathcal{F}_C(R)})$-exact by Lemma \ref{L2.14}. Then by the induction hypothesis, we can get
(2.5) is $\text{Hom}(-,\widetilde{\mathcal{F}_C(R)})$-exact. Dually, one gets another assertion. \endproof

The following result means that an iteration of the procedure used to define the $D_C$-projective
complexes yields exactly the $D_C$-projective complexes.

\begin{The}\label{T2.16} Let $X\in{\rm Ch}(R)$. Then $X$ is $D_C$-projective if and only if there exists a ${\rm Hom}(-,\widetilde{\mathcal{F}_C(R)})$-exact exact sequence of $D_{C}$-projective complexes
\vspace*{-2mm}$$\xymatrix{\mathbf{G}=\cdots \ar[r]&G_{1} \ar[r]^{\sigma_{1}}  &G_{0}  \ar[r]^{ \sigma_{0}} & G_{-1}  \ar[r]&\cdots}\vspace*{-2mm}$$
such that $X\cong$Coker$\sigma_{1}$.\end{The}

\proof $\Longrightarrow$) It is trivial.

$\Longleftarrow)$
Suppose that there exists a $\text{Hom}(-,\widetilde{\mathcal{F}_C(R)})$-exact exact sequence of $D_{C}$-projective complexes
\vspace*{-2mm}$$\xymatrix{\mathbf{G}=\cdots \ar[r]&G_{1} \ar[r]^{\sigma_{1}}&G_{0}  \ar[r]^{ \sigma_{0}} & G_{-1}  \ar[r]&\cdots}\vspace*{-2mm}$$
such that $X\cong{\rm Im}\sigma_0$.
Put $X_{i}={\rm Im}\sigma_i$ for any $i\in \mathbb{Z}$.  Then $X_0=X$ and we have ${\rm Hom}(-, \widetilde{\mathcal{F}_C(R)})$-exact exact sequence in ${\rm Ch}(R)$
\vspace*{-2mm}$$\xymatrix{0\ar[r]&X_{i+1}\ar[r]&G_i\ar[r]&X_i\ar[r]&0}\vspace*{-2mm}$$
for all $i\in \mathbb{Z}$. We wish to construct an exact sequence of complexes
satisfying the Definition \ref{D2.1}.

Consider the short exact sequence
\vspace*{-2mm}$$\xymatrix{0\ar[r]&X\ar[r]&G_{-1}\ar[r]&X_{-1}\ar[r]&0.}\vspace*{-2mm}$$
By Lemma \ref{L2.15}, there exist exact sequences
\vspace*{-2mm}$$\xymatrix{0\ar[r]&X\ar[r]&Q_{-1}\ar[r]&Y_{-1}\ar[r]&0}\vspace*{-2mm}$$ and
\vspace*{-2mm}$$\xymatrix{0\ar[r]&X_{-1}\ar[r]&Y_{-1}\ar[r]&V_{-1}\ar[r]&0}\vspace*{-2mm}$$
with $Q_{-1}\in \widetilde{\mathcal{P}_C(R)}$, $V_{-1}$  $D_C$-projective and the former one is ${\rm Hom}(-, \widetilde{\mathcal{F}_C(R)})$-exact.
Then by the pushout diagram
\vspace*{-2mm}$$\xymatrix{
 &                                        0 \ar[d]_{}                  & 0  \ar[d]_{}             &               &\\
 0 \ar[r]^{}        &X_{-1}\ar[d]^{} \ar[r]^{}     &Y_{-1}  \ar[d]_{}  \ar[r]^{ } & V_{-1} \ar@{=}[d]_{} \ar[r]^{ } &0  \\
 0 \ar[r]^{}        &G_{-2}\ar[d]^{ }\ar[r]^{}  &U_{-1} \ar[d]^{ } \ar[r]^{}   &V_{-1} \ar[r]^{}  &0  \\
 &               X_{-2} \ar[d]^{} \ar@{=}[r]        &X_{-2}  \ar[d]^{} &\\
                                       &  0               & 0   }\vspace*{-2mm}$$
we get an exact sequence
\vspace*{-2mm}$$\xymatrix{0\ar[r]&Y_{-1}\ar[r]&U_{-1}\ar[r]&X_{-2}\ar[r]&0.}\vspace*{-2mm}$$
By Corollary \ref{C2.10} and the exactness of the middle row in the above diagram, $U_{-1}$ is $D_C$-projective.
Since the first column in the above diagram is ${\rm Hom}(-, \widetilde{\mathcal{F}_C(R)})$-exact, ${\rm Ext}^1(X_{-2},\widetilde{\mathcal{F}_C(R)})=0$. It yields that $\xymatrix{0\ar[r]&Y_{-1}\ar[r]&U_{-1}\ar[r]&X_{-2}\ar[r]&0}$ is
${\rm Hom}(-, \widetilde{\mathcal{F}_C(R)})$-exact. By Lemma \ref{L2.15}, there exist exact sequences
\vspace*{-2mm}$$\xymatrix{0\ar[r]&Y_{-1}\ar[r]&Q_{-2}\ar[r]&Y_{-2}\ar[r]&0}\vspace*{-2mm}$$ and
\vspace*{-2mm}$$\xymatrix{0\ar[r]&X_{-2}\ar[r]&Y_{-2}\ar[r]&V_{-2}\ar[r]&0}\vspace*{-2mm}$$
with $Q_{-2}\in \widetilde{\mathcal{P}_C(R)}$, $V_{-2}$  $D_C$-projective and the former one is ${\rm Hom}(-, \widetilde{\mathcal{F}_C(R)})$-exact.
Then by the above argument, we have a ${\rm Hom}(-, \widetilde{\mathcal{F}_C(R)})$-exact exact sequence
\vspace*{-2mm}$$\xymatrix{0\ar[r]&Y_{-2}\ar[r]&U_{-2}\ar[r]&X_{-3}\ar[r]&0.}\vspace*{-2mm}$$
We proceed in this manner to get ${\rm Hom}(-, \widetilde{\mathcal{F}_C(R)})$-exact exact sequences
\vspace*{-2mm}$$\xymatrix{0\ar[r]&Y_{-i+1}\ar[r]&Q_{-i}\ar[r]&Y_{-i}\ar[r]&0}\vspace*{-2mm}$$
with $Q_{-i}\in \widetilde{\mathcal{P}_C(R)}$ for $i=1,2,\ldots$ where $Y_0=X$. Assembling these sequence together, we obtain a ${\rm Hom}(-, \widetilde{\mathcal{F}_C(R)})$-exact exact sequence
\vspace*{-2mm}$$\xymatrix{0\ar[r]& X\ar[r]& Q_{-1}\ar[r]& Q_{-2} \ar[r]&\cdots}         \eqno(\ast)\vspace*{-2mm}$$
with all $Q_{i}\in \widetilde{\mathcal{P}_C(R)}$.

Dually, we can get a ${\rm Hom}(-, \widetilde{\mathcal{F}_C(R)})$-exact exact sequence
\vspace*{-2mm}$$\xymatrix{\cdots\ar[r]& P_{1}\ar[r]& P_{0} \ar[r]&X\ar[r]&0}         \eqno(\ast\ast)\vspace*{-2mm}$$
with all $P_{i}\in \widetilde{\mathcal{P}(R)}$.

Finally, assembling the sequence $(\ast)$ and $(\ast\ast)$, we get a ${\rm Hom}(-, \widetilde{\mathcal{F}_C(R)})$-exact exact sequence
\vspace*{-2mm}$$\xymatrix{\mathbf{E}=\cdots\ar[r]^{} & P_{1} \ar[r]^{} & P_{0}\ar[r]^{} & Q_{-1}\ar[r]^{}& Q_{- 2}\ar[r]^{} &\cdots,}\vspace*{-2mm}$$
with all $P_{i}\in \widetilde{\mathcal{P}(R)}$ and all $Q_{i}\in \widetilde{\mathcal{P}_C(R)}$ such that $X\cong{\rm Im}(\xymatrix{P_{0} \ar[r]& Q_{-1}).}$  So $X$ is $D_C$-projective.\endproof

\begin{Cor}\label{2.17} Let $X\in{\rm Ch}(R)$. Then $X$ is $D_C$-projective if and only if there exists a ${\rm Hom}(-, \widetilde{\mathcal{F}_C(R)})$-exact exact sequence\vspace*{-2mm}$$\xymatrix{\cdots \ar[r]&W_{1} \ar[r]^{\sigma_{1}}  &W_{0}  \ar[r]^{ \sigma_{0}} & W_{-1}  \ar[r]&\cdots}\vspace*{-2mm}$$
in ${\rm Ch}(R)$ with all $W_i\in \widetilde{\mathcal{P}(R)}\cup\widetilde{\mathcal{P}_C(R)}$ such that $X\cong$Coker$\sigma_{1}$.\end{Cor}

\proof Immediate from Corollary \ref{C2.9} and Theorem \ref{T2.16}.\endproof

\section{$D_{C}$-Projective Dimension of Complexes}

Note that projective complexes are
$D_C$-projective by Corollary \ref{C2.9}, thus every complex admits a $D_{C}$-projective resolution. So we can define $D_{C}$-projective dimension of complexes as follows.

\begin{Def}\label{D3.1} Let $X\in {\rm{Ch}}(R)$. We will say that $X$ has $D_{C}$-projective dimension less than
or equal to $n$, denoted $D_{C}$-${\rm pd}(X)\leq n$, if there exists an exact sequence
$$\xymatrix{0\ar[r]&G_{n}\ar[r]&G_{n-1}\ar[r]&\cdots\ar[r]&G_{0}\ar[r]&X\ar[r]&0} $$ in ${\rm Ch}(R)$ with every $G_{i}$ being $D_{C}$-projective. If no such finite sequence exists, define $D_{C}$-${\rm pd}(X)=\infty$, otherwise, if $n$ is the least such integer, define $D_{C}$-${\rm pd}(X)= n$.\end{Def}

In this Section, we will give some criteria for computing $D_{C}$-${\rm pd}(X)$ of a complex $X$ if $D_{C}$-${\rm pd}(X)<\infty$. For this purpose, we need the following result.

\begin{Lem}\label{L3.2} Let $X\in {\rm Ch}(R)$. Consider two exact sequences
\vspace*{-2mm}
$$\xymatrix{ 0\ar[r]& H_{n}\ar[r]& G_{n-1} \ar[r]& G_{n-2}\ar[r]& \cdots\ar[r]&  G_{0}\ar[r]& X\ar[r]&  0 }\vspace*{-2mm}$$
and
$$\xymatrix{ 0\ar[r]& H_{n}'\ar[r]&G_{n-1}' \ar[r]&G_{n-2}'\ar[r]&\cdots\ar[r]&G_{0}'\ar[r]&X\ar[r]&0 }\vspace*{-2mm}$$
 with all $G_i,G_i'$ are $D_{C}$-projective complexes. Then $H_{n}$ is $D_{C}$-projective if and only if $H_{n}'$ is $D_{C}$-projective.\end{Lem}

\proof Using Corollaries \ref{C2.10} and \ref{C2.11}, the proof is similar to that of (i)$\Longrightarrow$(iii) in \cite[Theorem
1.2.7]{ch}.\endproof

A complex $X$ is said to have $\mathcal{P}_C$-projective dimension less than
or equal to $n$, denoted $\mathcal{P}_C$-${\rm pd}(X)\leq n$, if there is an exact sequence \vspace*{-2mm}$$\xymatrix{0\ar[r]& Q_{n}\ar[r]&Q_{n-1}\ar[r]&\cdots\ar[r]&Q_{1}\ar[r]&Q_{0}\ar[r]&X\ar[r]& 0}\vspace*{-2mm}$$ with each $Q_i\in \widetilde{\mathcal{P}_ C(R)}$ . If $n$ is the least then we set $\mathcal{P}_C$-${\rm pd}(X)=n$, and if there is no such $n$ then we set $\mathcal{P}_C$-${\rm pd}(X)=\infty$.
The $\mathcal{F}_C$-flat dimension of $X$, denoted by $\mathcal{F}_C$-${\rm pd}(X)$ can be
defined similarly.

\begin{The}\label{T3.3} Let $X\in{\rm Ch}(R)$ and $n\geq 0$. Then the following are equivalent:

(1) $D_{C}$-${\rm pd}(X)\leq n$.

(2) $D_{C}$-${\rm pd}(X)< \infty $ and ${\rm Ext}^m(X,H)$=0 for any $m>n$ and any $H\in {\rm Ch}(R)$ with $\mathcal{F}_C$-${\rm pd}(H)<\infty$.

(3) $D_{C}$-${\rm pd}(X)< \infty $ and ${\rm Ext}^m(X,H)$=0 for any $m>n$ and any $H\in \widetilde{\mathcal{F}_C(R)}$.

(4) For any exact sequence of complexes \vspace*{-2mm}$$\xymatrix{\cdots\ar[r]& G_{n}\ar[r]& G_{n-1} \ar[r]& \cdots\ar[r]& G_{0}\ar[r]& X\ar[r]&0 }\vspace*{-2mm}$$
with all $G_i$ $D_C$-projective, $K_{n}={\rm Ker}(\xymatrix{G_{n-1}\ar[r] &G_{n-2})}$ is $D_{C}$-projective.


(5) For any integer $t$ with $0\leq t\leq n$, there is an exact sequence of complexes\vspace*{-2mm}$$\xymatrix@-0.4pc {0\ar[r]& Q_{n}\ar[r]& \cdots\ar[r]& Q_{t+1}\ar[r]&  G_{t} \ar[r]& P_{t-1}\ar[r]&  \cdots\ar[r]& P_{0}\ar[r]& X\ar[r]&0}\vspace*{-2mm}$$such that $G_t$ is $D_C$-projective, $Q_i\in \widetilde{\mathcal{P}_C(R)}$ for $i> t$ and $P_i\in \widetilde{\mathcal{P}(R)}$
for $i< t$.
\end{The}

\proof (4)$\Longrightarrow$(1) and (2)$\Longrightarrow$(3) are trivial.

(1)$\Longrightarrow$(2) Since $D_{C}\text{-pd}(X)\leq n$, there is
an exact sequence
\vspace*{-2mm}$$\xymatrix{0 \ar[r]&G_{n}\ar[r]&G_{n-1} \ar[r]& \cdots\ar[r]&G_{0}\ar[r]&X\ar[r]&0 }\vspace*{-2mm}$$ with $G_i$ $D_C$-projective for all $0\leq i \leq n$. Then by dimension shifting and Lemma \ref{L2.2}, $\text{Ext}^m(X,H)\cong \text{Ext}^{m-n}(G_{n},H)$=0 for any $m>n$ and any $H\in {\rm Ch}(R)$ with $\mathcal{F}_C$-pd$(H)<\infty$.

(3)$\Longrightarrow$(4) Let
\vspace*{-2mm}$$\xymatrix{\cdots \ar[r]&G_{n}\ar[r]& G_{n-1} \ar[r]&\cdots\ar[r]&G_{0}\ar[r]&X\ar[r]&0}\vspace*{-2mm}$$ be an exact sequence 
in ${\rm Ch}(R)$ with all $G_i$ $D_{C}$-projective.
We will show that $K_{n}={\rm Ker}(\xymatrix{G_{n-1}\ar[r] &G_{n-2})}$ is $D_{C}$-projective.

By (3), we assume that $D_{C}\text{-pd}(X)=m<\infty $. Then there is
an exact sequence
\vspace*{-2mm}$$\xymatrix{0 \ar[r]&G_{m}'\ar[r]&G_{m-1}' \ar[r]&\cdots\ar[r]&G_{1}'\ar[r]&G_{0}'\ar[r]&X\ar[r]&0 }\vspace*{-2mm}$$ with all
$G_{i}'$ $D_{C}$-projective.
If $m \leq n$, there is nothing to prove. Now we assume that $m>n$. Set $K_{i}'={\rm Ker}(\xymatrix{G_{i-1}'\ar[r]& G_{i-2}'})$ for $i=1,2,\ldots m$, where $G_{-1}'=X$ and $K_{m}'=G_{m}'$.
By dimension shifting we have $\text{Ext}^{i+1}(X,H)\cong \text{Ext}^{1}(K_{i}',H)$ for any $i=n,n+1,\ldots,m-1$ and any $H\in \widetilde{\mathcal{F}_C(R)}$.
Thus $K_{n}'$ is $D_{C}$-projective inductively by (3) and Corollary \ref{C2.12}. Therefore $K_{n}$ is $D_{C}$-projective by Lemma \ref{L3.2}.

(1)$\Longrightarrow$(5) We proceed by induction on $n$.

If $n=1$, then there exists
an exact sequence $\xymatrix{0\ar[r]&D_1\ar[r] &D_0\ar[r]& X\ar[r]&0}$ with $D_0,D_1\,D_{C}$-projective.
By Lemma \ref{L2.14} with $A=0$, we get the exact sequences $\xymatrix{0\ar[r]&Q_1\ar[r] &G_0\ar[r]& X\ar[r]&0}$ and
$\xymatrix{0\ar[r]&G_1\ar[r] &P_0\ar[r]& X\ar[r]&0}$ with $G_0,G_1$ $D_{C}$-projective, $Q_1\in \widetilde{\mathcal{P}_C(R)}$ and $P_0\in \widetilde{\mathcal{P}(R)}$.

Next we suppose $n\geq2$. Then there exists an exact sequence of complexes
\vspace*{-2mm}$$\xymatrix{0\ar[r]& D_{n}\ar[r]&D_{n-1}\ar[r]&\cdots\ar[r]&D_{0}\ar[r]&X\ar[r]& 0}\eqno(\star)\vspace*{-2mm}$$ where all $D_i$ are $D_{C}$-projective.
Put $A={\rm Ker}\xymatrix{(D_1\ar[r]&D_0).}$ By applying Lemma \ref{L2.14} to the exact sequence
\vspace*{-2mm}$$\xymatrix{0\ar[r]&A\ar[r]& D_{1}\ar[r]&D_{0}\ar[r]&X\ar[r]& 0,}\vspace*{-2mm}$$ one gets the exactness of \vspace*{-2mm}$$\xymatrix{0\ar[r]&A\ar[r]& D_{1}'\ar[r]&P_{0}\ar[r]&X\ar[r]& 0}\vspace*{-2mm}$$ with $P_0\in \widetilde{\mathcal{P}(R)}$ and $D_1'$ $D_C$-projective. Hence we obtain the following exact sequence of complexes
\vspace*{-2mm}$$\xymatrix{0\ar[r]& D_{n}\ar[r]&D_{n-1}\ar[r]&\cdots\ar[r]&D_{2}\ar[r]&D_{1}'\ar[r]&P_{0}\ar[r]&X\ar[r]& 0.}\vspace*{-2mm}$$
Set $Y={\rm Ker}\xymatrix{(P_0\ar[r]&X).}$ Then $D_{C}\text{-pd}(Y)\leq n-1$. By the induction hypothesis, we can get an exact sequence
\vspace*{-2mm}$$\xymatrix@-0.3pc {0\ar[r]& Q_{n}\ar[r]& \cdots\ar[r]& Q_{t+1}\ar[r]&  G_{t} \ar[r]& P_{t-1}\ar[r]&  \cdots\ar[r]& P_{0}\ar[r]& X\ar[r]&0}\vspace*{-2mm}$$
with $G_t$ $D_{C}$-projective, all $P_i\in \widetilde{\mathcal{P}(R)}$ for $i<t$ and all $Q_i\in \widetilde{\mathcal{P}_C(R)}$ for $i>t$, where $1\leq t\leq n$.

Now it remains to show (5) for the case $t=0$. In the sequence $(\star)$, set $B={\rm Ker}\xymatrix{(D_0\ar[r]&X).}$ One gets the exactness of \vspace*{-2mm}$$\xymatrix{0\ar[r]& D_{n}\ar[r]&D_{n-1}\ar[r]&\cdots\ar[r]&D_{1}\ar[r]&B\ar[r]& 0.}\vspace*{-2mm}$$ By the induction
hypothesis, there is an exact sequence
\vspace*{-2mm}$$\xymatrix{0\ar[r]& Q_{n}\ar[r]& Q_{n-1}\ar[r]&\cdots\ar[r]&  Q_{2}\ar[r]& G_1'\ar[r]&B\ar[r]&0 ,}\vspace*{-2mm}$$
with $G_1'$ $D_{C}$-projective and all $Q_i\in \widetilde{\mathcal{P}_C(R)}$ with $2\leq i\leq n$. Set $A={\rm Coker}\xymatrix@-0.8pc{(Q_3\ar[r]&Q_2).}$ For the exact sequence
\vspace*{-2mm}$$\xymatrix{0\ar[r]&A\ar[r] &G_1'\ar[r]& D_0\ar[r]&X\ar[r]&0,}\vspace*{-2mm}$$ by Lemma \ref{L2.14}, we get an exact
sequence
\vspace*{-2mm}$$\xymatrix{0\ar[r]&A\ar[r] &Q_1\ar[r]& G_0\ar[r]&X\ar[r]&0}\vspace*{-2mm}$$ with $G_0$ $D_{C}$-projective and $Q_1\in \widetilde{\mathcal{P}_C(R)}$.
Thus we obtain the desired exact sequence
\vspace*{-2mm}$$\xymatrix{0\ar[r]& Q_{n}\ar[r]&Q_{n-1}\ar[r]&\cdots\ar[r]&Q_{2}\ar[r]&Q_{1}\ar[r]&G_{0}\ar[r]&X\ar[r]& 0}\vspace*{-2mm}$$
with $G_0$ $D_{C}$-projective and all $Q_i\in \widetilde{\mathcal{P}_C(R)}$ for $1\leq i\leq n$.

(5)$\Longrightarrow$(1) follows from Corollary \ref{C2.9}.\endproof

Using an argument as in the proof of \cite[Corollary 1.2.9]{ch}, we get the following
result by Lemma \ref{C2.10} and Theorem \ref{T3.3}.

\begin{Cor}\label{C3.4} Let $\xymatrix{0\ar[r]&X\ar[r] &Y\ar[r]& Z\ar[r]&0}$ be an exact sequence in ${\rm Ch}(R)$. Then the following hold:

(1) For any $n\geq 0$, if $D_{C}$-${\rm pd}(Z)\leq n$, then $D_{C}$-${\rm pd}(X)\leq n $ if and only if $D_{C}$-${\rm pd}(Y)\leq n$. Consequently,
$D_{C}$-${\rm pd}(X)\leq \max\{D_{C}$-${\rm pd}(Y),D_{C}$-${\rm pd}(Z)\}$
and
$D_{C}$-${\rm pd}(Y)\leq \max\{D_{C}$-${\rm pd}(X),D_{C}$-${\rm pd}(Z)\}.$

(2) If $D_{C}$-${\rm pd}(X)>D_{C}$-${\rm pd}(Z)$ or $D_{C}$-${\rm pd}(Y)>D_{C}$-${\rm pd}(Z),$  then
$D_{C}$-${\rm pd}(X)=D_{C}$-${\rm pd}(Y).$

(3) If $D_{C}$-${\rm pd}(Z)>0$ and $Y$ is $D_{C}$-projective, then $D_{C}$-${\rm pd}(X)=D_{C}$-${\rm pd}(Z)-1.$

In particular, if two complexes in the sequence $\xymatrix{0\ar[r]&X\ar[r] &Y\ar[r]& Z\ar[r]&0}$ have finite $D_{C}$-projective dimension, then
so is the third.\end{Cor}

Let $\mathcal{F}$ be a class of objects of an Abelian category $\mathcal{A}$ and $A$ an object of $\mathcal{A}$.
Following \cite{ejbook}, we say that a morphism $\xymatrix{f:F\ar[r]& A}$ is a $\mathcal{F}$-precover if $F\in\mathcal{F}$
and $\xymatrix{{\rm Hom}(F',F)\ar[r]&{\rm Hom}(F',A)\ar[r]&0}$
is exact for each $F'\in\mathcal{F}$. If such $f$ is an epimorphism, then we call $\xymatrix{f:F\ar[r]& A}$ is an epic $\mathcal{F}$-precover of $A$.  Recall that $A$ is said to have a special $\mathcal{F}$-precover if there is an exact
sequence $\xymatrix{0\ar[r]&K\ar[r] &F\ar[r]& A\ar[r]&0}$
with $F\in \mathcal{F}$ and ${\rm Ext}^1(\mathcal{F}, K)=0$. It is clear that $A$ has an epic $\mathcal{F}$-precover if it
has a special $\mathcal{F}$-precover. For more details about precovers, readers can refer to \cite{ejbook}.
The following result shows that a complex of finite $D_{C}$-projective
dimension can be approximated by a complex of finite C-projective dimension and
can also be approximated by a $D_{C}$-projective complex.

\begin{Cor}\label{C3.5} Let $X\in{\rm Ch}(R)$ with $D_{C}$-${\rm pd}(X)=n<\infty$. Then

(1) There exists an exact sequence $\xymatrix{0\ar[r]&X\ar[r] &Y\ar[r]& G\ar[r]&0}$ in ${\rm Ch}(R)$ with $G$ $D_{C}$-projective and $\mathcal{P}_C$-${\rm pd}(Y)=n$.

(2) $X$ admits a special $D_C$-projective precover $\xymatrix{0\ar[r]&K\ar[r] &G\ar[r]& X\ar[r]&0}$ with $\mathcal{P}_C$-${\rm pd}(K)=n-1$ if $n > 0$ and $K=0$ if $n=0$. \end{Cor}

\proof (1) If $X$ is $D_C$-projective then the result holds by Corollary \ref{C2.13}. Now
assume that $D_{C}$-${\rm pd}(X)=n>0$. Then we use Lemma \ref{L2.15}(1) with $A=0$ to get an exact sequence
$\xymatrix{0\ar[r]& X\ar[r]&Y\ar[r]&G\ar[r]&0}$ with $G$ $D_C$-projective and $\mathcal{P}_{C}$-${\rm pd}(Y)\leq n$.
By Corollary \ref{C3.4}(2), we have $D_{C}$-${\rm pd}(Y)=n$, and thus $\mathcal{P}_C$-${\rm pd}(Y)=n$.

(2) If $n=0$, it is trivial. Now assume that $n>0$. By Theorem \ref{T3.3}, there exists an exact sequence \vspace*{-2mm}$$\xymatrix{0\ar[r]& Q_{n}\ar[r]&Q_{n-1}\ar[r]&\cdots\ar[r]&Q_{2}\ar[r]&Q_{1}\ar[r]&G\ar[r]&X\ar[r]& 0}\vspace*{-2mm}$$
with $G$ $D_{C}$-projective and all $Q_i\in \widetilde{\mathcal{P}_C(R)}$ for $1\leq i\leq n$. Put $\xymatrix{K={\rm Ker}(G\ar[r]&X)}$. Then we have an exact sequence $\xymatrix{0\ar[r]&K\ar[r] &G\ar[r]& X\ar[r]&0}$ with $G$ $D_{C}$-projective and $\mathcal{P}_{C}$-${\rm pd}(K)\leq n-1$.
It follows from Corollary \ref{C3.4}(3) that $D_{C}$-${\rm pd}(K)=D_{C}$-${\rm pd}(X)-1=n-1$, and so $\mathcal{P}_C$-${\rm pd}(K)=n-1$. Also by Theorem \ref{T3.3}, ${\rm Ext}^1(G', K)=0$ for any $D_C$-projective complex $G'$. This completes the proof.
\endproof

\textbf{Acknowledgement.} The authors would like to appreciate the anonymous referee
for many helpful comments and suggestions.

{\small
{\em Authors' addresses}:
{\sc Yanhong Quan, Renyu Zhao, Chunxia Zhang}, Department of Mathematics, Northwest Normal University, Lanzhou,Gansu, China.
 E-mail: \texttt{786876393@qq.com, zhaory@nwnu.edu.cn, zhangcx@nwnu.edu.cn}.}

\end{document}